\documentclass{article}
\usepackage{amsthm}
\usepackage{amssymb}
\usepackage{amsmath}
\usepackage{amsfonts}
\usepackage{color}
\usepackage{indentfirst} %Ìí¼Ó
\usepackage{graphicx}
\usepackage[square, comma, sort&compress, numbers]{natbib}

\newtheorem{remark}{Remark}[section]

\newtheorem{lemma}{Lemma}[section]
\newtheorem{theorem}{Theorem}[section]

\newtheorem{corollary}{Corollary}[section]
\def\b1{\mbox{\boldmath $1$}}

\newenvironment{demo*}{\vspace{3mm}\noindent{\bf Proof.}}{\hfill $\Box$ \vspace{3mm}}

\begin{document}
\title{\bf \Large {Stochastic Orderings of Multivariate Elliptical Distributions  }}
{\color{red}{\author{%\normalsize{Ying Zhang}\\{\normalsize\it  School of Mathematical Sciences, Qufu Normal University}\\\noindent{\normalsize\it Shandong 273165, China}\\e-mail:  zhangying0513@gmail.com\\\and
\normalsize{Chuancun Yin\;\;  }\\
%\thanks{Corresponding author.}\\
{\normalsize\it  School of Statistics,  Qufu Normal University}\\
\noindent{\normalsize\it Shandong 273165, China}\\
e-mail:  ccyin@qfnu.edu.cn \\
}}}%\\
\maketitle
\vskip0.01cm
%\newpage
{\large {\bf Abstract}}  Let ${\bf X}$  and ${\bf Y}$ be two $n$-dimensional elliptical random vectors, we establish an identity for
$E[f({\bf Y})]-E[f({\bf X})]$, where $f: \Bbb{R}^n \rightarrow \Bbb{R}$   satisfies some regularity conditions. Using this identity we provide a unified derivation of sufficient and necessary conditions for classifying multivariate elliptical random vectors according to several main integral stochastic orders. As a consequence we obtain new
inequalities by applying it to multivariate elliptical distributions. The results generalize the corresponding ones for multivariate normal random vectors in the literature.

\medskip



%\newpage
%\noindent{\bf 1.~~Introduction}
\numberwithin{equation}{section}
\section{Introduction}\label{intro}
Stochastic orders   provide methods of comparing random variables and vectors which are now used in many areas such as  statistics and probability (Cal and Carcamo (2006), Hu and Zhuang (2006),   M\"uller and   Scarsini (2006),  Fill and  Kahn (2013) and  Hazra et al. (2017)),  operations research (F\'abi\'an et al. (2011)),  actuarial sciences and  economic theory (Ba\"uerle and  Bayraktar (2014),  L\'opez-D\'iaz et al. (2018)),   risk management and other related fields  (B\"auerle  and M\"uller (2006)). For a comprehensive review of the properties and characterizations of stochastic orderings, including a variety of applications, the reader is referred to the monographs of M\"uller and Stoyan (2002), Denuit et al. (2005),  and Shaked and Shanthikumar (2007).
Many of these orders  are characterized by the so-called integral stochastic
orders which is obtained by comparing expectations of functions in a certain class.  A general   treatment
for these orders  has been given in Whitt (1986) and M\"uller (1997).

 Elliptical distributions are generalizations of the multivariate normal distributions and, therefore, share many of the tractable properties. This class of distributions was introduced by Kelker (1970) and was extensively discussed in Fang et al. (1987).  This generalization of the normal family seems to provide an attractive tool for  statistics, economics and finance,     which can describe fat or
light tails and tail dependence among components of a vector.  Interested
readers are referred to monograph  of Gupta, Varga and Bodnar (2013) and some recent papers  of El Karoui (2009),  Hu et al. (2019) and Sha et al. (2019).   M\"uller (2001) studied  the stochastic ordering characterizations of multivariate normal random vectors.   Arlotto and Scarsini (2009) unified and generalized several known results on comparisons of multivariate normal random vectors in the sense of different stochastic orders by introducing the so-called Hessian order.
  Landsman and Tsanakas (2006) derived necessary and sufficient conditions for classifying bivariate elliptical distributions through the concordance ordering. Ding and Zhang (2004) extended  the results in   M\"uller (2001)  to Kotz-type distributions   which form an important specially  class of elliptical symmetric distributions.
 Necessary and sufficient conditions for  convex order and increasing convex order of
general  multivariate elliptical random vectors  had not been found until the work of  Pan et al. (2016). However, few results can be found in the literature that characterize the supermodular order of multivariate elliptical distributions. It is the aim of this paper to fill this gap.  We will give some sufficient and necessary conditions for  supermodular order  of multivariate elliptical random vectors.  For the known results  such as on the convex ordering and the increasing convex ordering of multivariate elliptical random vectors, we provide a different simple proof.

 The rest of the paper is organized as follows. Section 2   recalls  some useful notions
that will be used in the sequel, such as certain properties of  stochastic orders and elliptical distributions.
 Section 3   presents   necessary  and sufficient conditions  for several important  stochastic orders   of  multivariate
elliptical distributions. Section 4 provides    two applications of the  main results.

\section{Preliminaries}

Throughout this paper we  use the following notations.
We use bold letters to denote vectors or matrices. For example, ${\bf X}' =(X_1,\cdots, X_n)$ is a row vector and
${\bf \Sigma} = (\sigma_{ij})_{n\times n}$  is an $n\times n$ matrix. In particular, the symbol ${\bf 0}_n$ denotes
the  $n$-dimensional column vector with all entries equal to 0,  ${\bf 1}_{n}$  denotes the $n$-dimensional column
vector with all components equal to  $1$,  and ${\bf 1}_{n\times n}$ denotes the ${n\times n}$ matrix with all entries equal to $1$.   Denote by ${\bf O}_{n\times n}$ the $n \times n$ matrix  with all entries  equal to $0$  and  by ${\bf I}_n$  the $n\times n$ identity matrix. For symmetric matrices $A$ and $B$ of the same size, the notion $A\preceq B$ or $B-A  \succeq {\bf O}$ means that $B-A$ is positive semi-definite. The inequality between vectors or matrices denotes componentwise inequalities.
Throughout this paper, the terms of increasing and decreasing are used in the weak sense.  All integrals and expectations are implicitly assumed to exist whenever they appear.

\subsection{Stochastic orders}
 In this section, we summarize some important definitions and facts about the
stochastic orderings of random vectors.
 The standard references for stochastic orderings are the monographs by
 Denuit et al.   (2005) and  Shaked and Shanthikumar (2007).
 For a function  $f: \Bbb{R}^n \rightarrow \Bbb{R}$,  ${\bf x}\in \Bbb{R}^n$,  $i\in \{1,\cdots,n\}$ and $\delta> 0$, we define
 the difference operator $\Delta_i^{\delta}$ as
 $$\Delta_i^{\delta}f({\bf x}) = f({\bf x}+\delta {\bf e}_i)-f ({\bf x}),$$
 where ${\bf e}_i=(0,\dots, 0, 1, 0,\dots, 0)$ denotes the $i$th unit vector.
 In case $n = 1$ we simply write
  $$\Delta^{\delta}f( x) = f(x+\delta)-f(x).$$
 A function  $f: \Bbb{R}^n \rightarrow \Bbb{R}$     is said to be increasing, if
 $\Delta_i^{\delta}f({\bf x})\ge 0$, for all ${\bf x}\in \Bbb{R}^n$, $\delta>0$ and $i =1,\cdots,n$.
 A function  $f: \Bbb{R}^n \rightarrow \Bbb{R}$   is  is said to be   supermodular, if
 $\Delta_i^{\delta} \Delta_j^{\varepsilon}f({\bf x})\ge 0$, for all ${\bf x}\in \Bbb{R}^n$, $\delta, \varepsilon>0$ and  $1\le i<j\le n$.
Equivalently, a function $f:  {\Bbb{R}}^n \rightarrow {\Bbb R}$ is said to be supermodular
if for any ${\bf x, y}\in {\Bbb R}^n$ it holds that
$$f({\bf x}) + f({\bf y}) \le f ({\bf x} \wedge {\bf y}) + f ({\bf x} \vee {\bf y}),$$
where the operators $\wedge$ and  $\vee$ denote coordinatewise minimum and
maximum, respectively.  A function $f$ is supermodular if and only if $-f$
is submodular.
A function  $f: \Bbb{R}^n \rightarrow \Bbb{R}$  is said to be componentwise convex if $f$ is convex in each
argument when the other arguments are hold fixed.
A function  $f: \Bbb{R}^n \rightarrow \Bbb{R}$   is   said to be   directionally convex, if
 $\Delta_i^{\delta} \Delta_j^{\varepsilon}f({\bf x})\ge 0$, for all ${\bf x}\in \Bbb{R}^n$, $\delta, \varepsilon>0$ and  $1\le i, j\le n$.
That is  $f: \Bbb{R}^n \rightarrow \Bbb{R}$   is directionally convex if it is supermodular and componentwise convex.
Directional convexity neither implies, nor is implied by, conventional convexity.
The supermodular order compares only the dependence structure of vectors
with fixed equal marginals, whereas the increasing directionally convex order also compares the
marginals both invariability and location, where the marginals are possibly different.
 However,
a univariate function is directionally convex if, and only if, it is convex.
A function  $f: \Bbb{R}^n \rightarrow \Bbb{R}$   is   said to be   $\Delta$-monotone function, if for all
$\{i_1, \cdots,i_k\}\subset \{1, \cdots, n\}$ and every $\delta_1,\cdots,\delta_k>0$,
 $\Delta_{i_1}^{\delta_1}\cdots \Delta_{i_k}^{\delta_{k}}f({\bf x})\ge 0$, for all ${\bf x}\in \Bbb{R}^n$, $\delta, \varepsilon>0$.

Let us now recall the definitions of stochastic orders that will be used later.

Let ${\cal{F}}$ be some class of measurable functions $f:  {\Bbb{R}}^n \rightarrow {\Bbb R}$, for two random vectors ${\bf X}$ and ${\bf Y}$ in  ${\Bbb{R}}^n$, we say that ${\bf X}\le_{\cal F}{\bf Y}$ if $E[f({\bf X})]\le E[f({\bf Y})]$ holds for all $f\in {\cal{F}}$ whenever the expectation
  is well defined.  We list   a few important examples as follows.

$\bullet$ Usual stochastic order:  ${\bf X}\le_{st} {\bf Y}$,  if $E[f({\bf X})]\le E[f({\bf Y})]$  for all increasing
functions $f: \Bbb{R}^n \rightarrow \Bbb{R}$.

$\bullet$ Convex order:  ${\bf X}\le_{cx} {\bf Y}$,  if $E[f({\bf X})]\le E[f({\bf Y})]$  for all  convex
functions $f: \Bbb{R}^n \rightarrow \Bbb{R}$.

$\bullet$  Linear convex order: ${\bf X}\le_{lcx} {\bf Y}$,   if $E[f({\bf a'X})]\le E[f({\bf a'Y})]$  for all ${\bf a}\in\Bbb{R}^n$ and for all  convex functions $f: \Bbb{R}^n \rightarrow \Bbb{R}$.

$\bullet$ Increasing convex order:  ${\bf X}\le_{icx} {\bf Y}$,  if $E[f({\bf X})]\le E[f({\bf Y})]$  for all increasing convex
functions $f: \Bbb{R}^n \rightarrow \Bbb{R}$.

$\bullet$ Componentwise convex order: ${\bf X}\le _{ccx} {\bf Y}$, if $E[f({\bf X})]\le E[f({\bf Y})]$ for all componentwise convex functions $f: \Bbb{R}^n \rightarrow \Bbb{R}$.

$\bullet$  Increasing componentwise convex order: ${\bf X}\le _{iccx} {\bf Y}$, if $E[f({\bf X})]\le E[f({\bf Y})]$ for all  increasing componentwise convex functions $f: \Bbb{R}^n \rightarrow \Bbb{R}$.

$\bullet$  Supermodular order: ${\bf X}\le_{sm} {\bf Y}$,  if $E[f({\bf X})]\le E[f({\bf Y})]$ for all   supermodular functions $f: \Bbb{R}^n \rightarrow \Bbb{R}$.

 $\bullet$  Increasing supermodular order: ${\bf X}\le_{ism} {\bf Y}$,  if $E[f({\bf X})]\le E[f({\bf Y})]$ for all  increasing supermodular functions $f: \Bbb{R}^n \rightarrow \Bbb{R}$.

$\bullet$    Directionally convex order:  ${\bf X}\le_{dcx} {\bf Y}$, if    $E[f({\bf X})]\le E[f({\bf Y})]$  for all   directionally convex functions $f: \Bbb{R}^n \rightarrow \Bbb{R}$.

$\bullet$   Increasing directionally convex:  ${\bf X}\le_{idcx} {\bf Y}$, if    $E[f({\bf X})]\le E[f({\bf Y})]$  for all increasing  directionally convex functions $f: \Bbb{R}^n \rightarrow \Bbb{R}$.

%$\bullet$ Upper orthant order:  ${\bf X}\le_{uo} {\bf Y}$,   if $E[f({\bf X})]\le E[f({\bf Y})]$  for all    $\Delta$-monotone  functions $f: \Bbb{R}^n \rightarrow \Bbb{R}$.

%$\bullet$ Componentwise convex order: ${\bf X}\le _{ccx} {\bf Y}$,  if $E[f({\bf X})]\le E[f({\bf Y})]$  for all componentwise convex functions $f: \Bbb{R}^n \rightarrow \Bbb{R}$.

%$\bullet$ Increasing componentwise convex order: ${\bf X}\le _{iccx} {\bf Y}$,  if $E[f({\bf X})]\le E[f({\bf Y})]$  for all increasing componentwise convex functions $f: \Bbb{R}^n \rightarrow \Bbb{R}$.

The componentwise convex order   was introduced in Mosler (1982),  the  directionally convex   was introduced in  Shaked and  Shanthikumar (1990) and
 the increasing directionally convex   was introduced in  Meester and  Shanthikumar (1993).

 For a random vector ${\bf X}=(X_1, \cdots, X_n)$, we denote by
$$F_{\bf X} ({\bf t}):=P({\bf X}\le {\bf t})=P(X_1\le t_1,\cdots, X_n\le t_n), {\bf t}=(t_1,\cdots, t_n)\in {\Bbb R}^n,$$
and
$${\overline F}_{\bf X} ({\bf t}):=P({\bf X}> {\bf t})=P(X_1> t_1,\cdots, X_n> t_n), {\bf t}=(t_1,\cdots, t_n)\in {\Bbb R}^n,$$
the multivariate distribution function and the multivariate survival function, respectively.
The following definition is taking from M\"uller and Scarsini (2000).

{\bf Definition 2.4.} (a) Assume that  ${\bf X, Y}\in {\Bbb R}^n$ are two random vectors. \\
(a) \; ${\bf X}$ is said to be smaller than   ${\bf Y}$ in the upper orthant order,
written ${\bf X}\le _{uo} {\bf Y}$, if  ${\overline F}_{\bf X} ({\bf t})\le {\overline F}_{\bf Y} ({\bf t})$ for all ${\bf t}\in {\Bbb R}^n$.\\
(b) \;  ${\bf X}$ is said to be smaller than   ${\bf Y}$ in the upper orthant order,
written ${\bf X}\le _{lo} {\bf Y}$, if  ${F}_{\bf X} ({\bf t})\le {F}_{\bf Y} ({\bf t})$ for all ${\bf t}\in {\Bbb R}^n$.\\
(c) \; ${\bf X}$ is said to be smaller than   ${\bf Y}$ in the concordance order,
 written ${\bf X}\le _{c} {\bf Y}$, if both   ${\bf X}\le _{uo} {\bf Y}$ and  ${\bf X}\le _{lo} {\bf Y}$ hold.

 The orthant orders have been treated by Shaked and Shanthikumar
(1994) and the concordance order was introduced by Joe (1990).   We have the implication ${\bf X}\le_{sm} {\bf Y}$   $\Rightarrow$ ${\bf X}\le _{uo} {\bf Y}$ and ${\bf X}\le _{lo} {\bf Y}$, and hence also  ${\bf X}\le_{sm} {\bf Y}$   $\Rightarrow$ ${\bf X}\le _{c} {\bf Y}$.

The upper orthant order can be defined alternatively by $\Delta$-monotone functions. The following lemma can be founded in  R\"uschendorf (1980).
\begin{lemma}
 ${\bf X}\le _{uo} {\bf Y}$ if and only if   $E[f({\bf X})]\le E[f({\bf Y})]$ holds for all $\Delta$-monotone functions $f: \Bbb{R}^n \rightarrow \Bbb{R}$.
\end{lemma}
  The following necessary and sufficient conditions for several important stochastic orders can be found in Denuit and M\"uller (2002) and Arlotto and Scarsini (2009).

1.  ${\bf X}\le_{sm} {\bf Y}$  if, and only if $E[f({\bf X})]\le E[f({\bf Y})]$ holds for all twice differentiable
functions $f: \Bbb{R}^n \rightarrow \Bbb{R}$ satisfying $\frac{\partial^2 }{\partial x_i \partial x_j}f({\bf x})\ge 0$ for ${\bf x}\in \Bbb{R}^n$ and all $1\le i<j\le n$.

2.  ${\bf X}\le_{ism} {\bf Y}$  if, and only if $E[f({\bf X})]\le E[f({\bf Y})]$  holds for all twice differentiable
functions $f: \Bbb{R}^n \rightarrow \Bbb{R}$ satisfying $\frac{\partial }{\partial x_i}f({\bf x})\ge 0$ for ${\bf x}\in \Bbb{R}^n$ and all $1\le i\le n$, and $\frac{\partial^2 }{\partial x_i \partial x_j}f({\bf x})\ge 0$ for ${\bf x}\in \Bbb{R}^n$ and all $1\le i<j\le n$.

3. ${\bf X}\le_{dcx} {\bf Y}$ if, and only if,    $E[f({\bf X})]\le E[f({\bf Y})]$  holds for all twice differentiable
functions $f: \Bbb{R}^n \rightarrow \Bbb{R}$ satisfying    $\frac{\partial^2 }{\partial x_i \partial x_j}f({\bf x})\ge 0$ for ${\bf x}\in \Bbb{R}^n$ and all $1\le i,j\le n$.

4. ${\bf X}\le_{idcx} {\bf Y}$  if, and only if $E[f({\bf X})]\le E[f({\bf Y})]$  holds for all twice differentiable
functions $f: \Bbb{R}^n \rightarrow \Bbb{R}$ satisfying $\frac{\partial }{\partial x_i}f({\bf x})\ge 0$ for ${\bf x}\in \Bbb{R}^n$ and all $1\le i\le n$, and $\frac{\partial^2 }{\partial x_i \partial x_j}f({\bf x})\ge 0$ for ${\bf x}\in \Bbb{R}^n$ and all $1\le i, j\le n$.

5. ${\bf X}\le_{uo} {\bf Y}$  if, and only if $E[f({\bf X})]\le E[f({\bf Y})]$ holds for all infinitely  differentiable
functions $f: \Bbb{R}^n \rightarrow \Bbb{R}$ satisfying $\frac{\partial^k }{\partial x_{i_1}\cdots \partial x_{i_k}}f({\bf x})\ge 0$ for ${\bf x}\in \Bbb{R}^n$ and all $1\le {i_1}<\cdots\le i_k\le n$.

6. ${\bf X}\le_{ccx} {\bf Y}$  if, and only if $E[f({\bf X})]\le E[f({\bf Y})]$  holds for all twice differentiable
functions $f: \Bbb{R}^n \rightarrow \Bbb{R}$ satisfying $\frac{\partial^2 }{\partial x_i^2}f({\bf x})\ge 0$ for ${\bf x}\in \Bbb{R}^n$ and all $1\le i\le n$.

7. ${\bf X}\le_{iccx} {\bf Y}$  if, and only if $E[f({\bf X})]\le E[f({\bf Y})]$  holds for all twice differentiable
functions $f: \Bbb{R}^n \rightarrow \Bbb{R}$ satisfying  $\frac{\partial }{\partial x_i}f({\bf x})\ge 0$   and  $\frac{\partial^2 }{\partial x_i^2}f({\bf x})\ge 0$ for ${\bf x}\in \Bbb{R}^n$ and all $1\le i\le n$.

We first list the results of stochastic orderings for univariate elliptical distributions.
For the case of univariate normal distributions can be found in M\"uller (2001).
\begin{lemma}  Let $X\sim E_1(\mu_x, \sigma_y^2, \phi)$ and $Y\sim E_1(\mu_y, \sigma_y^2, \phi)$. Then

(i)\; $X \le_{st} Y$ if and only if $\mu_x\le \mu_y$ and $\sigma_x=\sigma_y$, provided that  $X$ and $Y$ supported on $\Bbb{R}$  ( Davidov and Peddada(2013));

(ii)\; $X \le_{cx} Y$ if and only if $\mu_x=\mu_y$ and $\sigma_x\le \sigma_y$ (Pan et al. (2016));

(iii)\; $X \le_{icx} Y$ if and only if $\mu_x\le \mu_y$ and $\sigma_x\le \sigma_y$ (Pan et al. (2016)).
\end{lemma}

Now we list the results of stochastic orderings for   multivariate elliptical distributions.
\begin{lemma} (Pan et al. (2016)) Let  ${\bf{X}}\sim E_n ({\boldsymbol \mu}_x,{\bf \Sigma}_x,\phi)$ and   ${\bf{Y}}\sim E_n ({\boldsymbol \mu}_y,{\bf \Sigma}_y,\phi)$.  Then the following statements are equivalent:

(1) ${\boldsymbol \mu}_x={\boldsymbol \mu}_y$
and ${\bf \Sigma}_y-{\bf \Sigma}_x$ is positively semi-definite;

(2) ${\bf X}\le_{lcx} {\bf Y}$;

(3) ${\bf X}\le_{cx} {\bf Y}$.
\end{lemma}
For the case of increasing convex order, the sufficient and necessary conditions  seem to be unknown. The following  sufficient condition for the increasing convex order can be found in Pan et al. (2016).  The results for the case of multivariate normal distributions  can be found in M\"uller (2001).
\begin{lemma}  Let  ${\bf{X}}\sim E_n ({\boldsymbol \mu}_x,{\bf \Sigma}_x,\phi)$ and   ${\bf{Y}}\sim E_n ({\boldsymbol \mu}_y,{\bf \Sigma}_y,\phi)$.

(i) If  ${\boldsymbol \mu}_x\le {\boldsymbol \mu}_y$ and ${\bf \Sigma}_y-{\bf \Sigma}_x$ is positively semi-definite, then
${\bf X}\le_{icx} {\bf Y}$.

(ii) If ${\bf X}\le_{icx} {\bf Y}$, then ${\boldsymbol \mu}_x\le {\boldsymbol \mu}_y$ and ${\bf a}'({\bf \Sigma}_y-{\bf \Sigma}_x){\bf a}\ge 0$ for all ${\bf a}\ge 0.$
\end{lemma}

\subsection{  Some background on the elliptical distributions}
 The class of multivariate elliptical distributions is a natural extension to the class of
multivariate Normal distributions.
 We follow the notation of  Cambanis et al.  (1981) and  Fang et al. (1990).
An $n \times 1$ random vector $X = (X_1, X_2,\cdots, X_n)'$ is said to have an elliptically
symmetric distribution if its characteristic function has the form $e^{i{\bf t}'{\boldsymbol \mu}}\phi({\bf t}'{\bf \Sigma}{\bf t})$ for all ${\bf t}\in \Bbb{R}^n $,  denoted ${\bf{X}}\sim E_n ({\boldsymbol \mu},{\bf \Sigma},\phi)$,
  where  $\phi  \in {\bf \Psi}_n$ is called the characteristic generator satisfying $\phi(0)=1$,
$\boldsymbol{\mu}$ ($n$-dimensional vector) is its location parameter    and  $\bf{\Sigma}$ ($n\times n$ matrix with $\bf{\Sigma}\succeq {\bf O}$) is its dispersion matrix (or scale matrix). The mean vector $E({\bf X})$   (if it exists)
coincides with the location vector and the covariance matrix  Cov$({\bf X})$ (if it exists), being $-2\phi'(0){\bf \Sigma}$.
It is interesting to note that in the one-dimensional case, the class of elliptical distributions consists mainly of the class of symmetric distributions which include well-known distributions like Normal and Student $t$  distributions.
 It is well known that
$\bf X$ admits the stochastic representation
\begin{equation}
{\bf X}={\boldsymbol \mu}+R{\bf A}'{\bf U}^{(n)},
\end{equation}
where ${\bf A}$  is a square matrix such that ${\bf A}'{\bf A}= {\bf \Sigma}$, ${\bf U}^{(n)}$  is uniformly distributed on the unit sphere ${\cal S}^{n-1}=\{{\bf u}\in \Bbb{R}^n: {\bf u}'{\bf u}=1\} $,  $R\ge 0$   is the random variable with $R \sim F$ in $[0, \infty)$ called the generating variate and $F$ is called the generating distribution function, $R$ and  ${\bf U}^{(n)}$ are  independent. The mean vector $E({\bf X})$ exists if and only if $E(R)$ exists and $E({\bf X})={\boldsymbol \mu}$;  The covariance matrix  Cov$({\bf X})$  exists if and only if $E(R^2)$ exists and Cov$({\bf X})=\frac{1}{n} E(R^2){\bf \Sigma}.$
In general,   an elliptically distributed random vector ${\bf{X}}\sim E_n ({\boldsymbol \mu},{\bf \Sigma},\phi)$ does not necessarily
possess a density. It is well known that ${\bf{X}}$  has a density if and only if $R$ has
a density and ${\bf \Sigma}\succ {\bf O}$. The density has the form
\begin{equation}
f({\bf x})=c_n|{\bf \Sigma}|^{-\frac{1}{2}}g_n(({\bf x}- {\boldsymbol \mu})'{\bf \Sigma}^{-1}({\bf x-{\boldsymbol \mu} })), \;{\bf x}\in \Bbb{R}^n,
\end{equation}
for some nonnegative function $g_n$  called the density generator and for some constant
$c_n$ called the normalizing constant. One sometimes
writes $X \sim E_n ({\boldsymbol \mu},{\bf \Sigma},g_n)$ for the $n$-dimensional elliptical distributions generated from the
function $g_n$.

The class of elliptical distributions possesses the linearity property. Consider the affine transformations of the form ${\bf Y = BX+b}$ of a random vector  ${\bf{X}}\sim E_n({\boldsymbol \mu},{\bf \Sigma},\phi)$, where  $\bf B$ is a $m\times n$ matrix with $m<n$ and $rank({\bf B})=m$ and ${\bf b}\in\Bbb{R}^m$. Then  ${\bf Y} \sim E_n ({\bf B}{\boldsymbol \mu}+{\bf b}, {\bf B}{\bf \Sigma}{\bf B}',\phi)$.
 Taking $B=(\alpha_1,\cdots,\alpha_n):={\boldsymbol \alpha}'$  leads to
$${\boldsymbol \alpha}'{\bf X}\sim E_1({\boldsymbol \alpha}'{\boldsymbol \mu},  {\boldsymbol \alpha}'{\bf \Sigma}{\boldsymbol \alpha},\phi).$$
In particular,
$$X_k\sim Ell_1(\mu_k,\sigma_k^2, \phi), \;\; k=1,2,\cdots, n,$$
 and
 $$\sum_{k=1}^n X_k\sim Ell_1\left(\sum_{k=1}^n\mu_k,\sum_{k=1}^n\sum_{l=1}^n \sigma{_{kl}}, \phi\right).$$

\subsection{An identity for multivariate elliptical distributions}

If $f: \Bbb{R}^n \rightarrow \Bbb{R}$ is twice continuously differentiable, we write as usual
$$\nabla f({\bf x})=\left(\frac{\partial }{\partial x_1}f({\bf x)},\cdots, \frac{\partial }{\partial x_n}f({\bf x)}\right)',\;
 {\bf H}_f({\bf x})=\left(\frac{\partial^2}{\partial x_i \partial x_i}f({\bf x)} \right)_{n\times n}$$
  for the gradient and the Hessian matrix of $f$. It is well known that  $f$ is convex if and only if  ${\bf H}_f({\bf x})$
is positive semidefinite for any ${\bf x}\in \Bbb{R}^n$;   $f$ is strictly convex if and only if  ${\bf H}_f({\bf x})$
is positive definite for any ${\bf x}\in \Bbb{R}^n$. A function is supermodular if and
only if its Hessian has nonnegative off-diagonal elements, i.e.
  $f$ is supermodular if and only if $\frac{\partial^2}{\partial x_i \partial x_i}f({\bf x)}\ge 0$ for every $i\neq j$ and ${\bf x}\in {\Bbb{R}}^n$ (c.f.  Carter (2001, Proposition 4.2)).

The following  two results  can be found in    Houdr\'e et al. (1998),  M\"uller (2001) and Denuit and M\"uller  (2002) in the multivariate normal case.  Ding and Zhang (2004) extended the result from
multivariate normal distributions to Kotz-type distributions which form an important class of elliptically symmetric
distributions.   We develop an identity for multivariate elliptical distributions   which may be of independent interest.

\begin{lemma} Let ${\bf{X}}\sim Ell_n ({\boldsymbol \mu}^x,{\bf \Sigma}^{x},\psi)$ and ${\bf{Y}}\sim Ell_n ({\boldsymbol \mu}^y,{\bf \Sigma}^{y}, \psi)$, with
${\bf \Sigma}^{x}$ and ${\bf \Sigma}^{y}$ positive definite. Let $\phi_{\lambda}$ be the density
function of
$$
Ell_n (\lambda{\boldsymbol \mu}^y+(1-\lambda){\boldsymbol \mu}^x , \lambda{\bf \Sigma}^{y}+(1-\lambda){\bf \Sigma}^{x}, \psi),\; 0\le \lambda \le 1,             $$
and $\phi_{1\lambda}$ be the density
function of
$$
Ell_n (\lambda{\boldsymbol \mu}^y+(1-\lambda){\boldsymbol \mu}^x , \lambda{\bf \Sigma}^{y}+(1-\lambda){\bf \Sigma}^{x}, \psi_1),\; 0\le \lambda \le 1,             $$
where
$$\psi_1(u)=\frac{1}{E(R^2)}\int_0^{\infty} {}_0F_1\left(\frac{n}{2}+1;-\frac{r^2 u}{4}\right)r^2P(R\in dr).$$
Here
$${}_0F_1(\gamma;z)=\sum_{k=0}^{\infty}\frac{1}{(\gamma)_k}\frac{z^k}{k!},$$
is the  generalized hypergeometric series of order $(0, 1)$, $R$ is defined by (2.1) with $E(R^2)<\infty$.
Moreover, assume that
$f: \Bbb{R}^n \rightarrow \Bbb{R}$ is twice continuously differentiable and satisfies some polynomial
growth conditions at infinity:
$$f({\bf x})=O(||{\bf x}||),\; \bigtriangledown f({\bf x})=O(||{\bf x}||).$$
 Then
\begin{eqnarray*}
 E[f({\bf Y})]-E[f({\bf X})]&=&\int_0^1\int_{\Bbb{R}^n}({\boldsymbol \mu}^y-{\boldsymbol \mu}^x)'\nabla f({\bf x})\phi_{\lambda}({\bf x})d{\bf x}d\lambda\\
&&+\frac{E(R^2)}{2n}\int_0^1\int_{\Bbb{R}^n}tr\{({\bf \Sigma}^{y}-{\bf \Sigma}^{x}){\bf H}_f({\bf x})\}\phi_{1\lambda}({\bf x})d{\bf x}d\lambda,
\end{eqnarray*}
where tr$(A)$ denotes the trace of the matrix $A$.
\end{lemma}

{\bf Proof} See Appendix 2.

Using Lemma 2.5 and the same argument as  in the  proof of Corollary 3 in  M\"uller (2001) we have

 \begin{corollary} Let ${\bf{X}}\sim Ell_n ({\boldsymbol \mu}^x,{\bf \Sigma}^{x},\psi)$ and ${\bf{Y}}\sim Ell_n ({\boldsymbol \mu}^y,{\bf \Sigma}^{y},\psi)$, with
${\bf \Sigma}^{x}$ and ${\bf \Sigma}^{y}$ positive definite or  positive semidefinite, and   assume that
$f: \Bbb{R}^n \rightarrow \Bbb{R}$   satisfies the conditions of Lemma 2.5.
Then $E[f({\bf X})]\le E[f({\bf Y})]$ if the following two conditions hold for all ${\bf x}\in \Bbb{R}^n$:
$$\sum_{i=1}^n({\mu}^{y}_{i}-{ \mu}^{x}_{i})\frac{\partial }{\partial x_i}f({\bf x})\ge 0,$$
and
$$\sum_{i=1}^n\sum_{j=1}^n({\sigma}^y_{ij}-{\sigma}^x_{i j})\frac{\partial^2 }{\partial x_i \partial x_j}f({\bf x})\ge 0.$$
\end{corollary}

\section{Main results}

The following results can be extracted from Davidov and Peddada (2013). The   multivariate normal  case   can be found in M\"uller (2001). Here we provide a different proof.
\begin{theorem}   Let  ${\bf{X}}\sim E_n ({\boldsymbol \mu}^x,{\bf \Sigma}^x,\phi)$ and   ${\bf{Y}}\sim E_n ({\boldsymbol \mu}^y,{\bf \Sigma}^y,\phi)$  be   two $n$-dimensional elliptically
distributed random vectors supported on ${\Bbb R}^n$. Then ${\bf X}\le_{st} {\bf Y}$ if and only if  ${\boldsymbol \mu}^x\le {\boldsymbol \mu}^y$ and ${\bf \Sigma}^y={\bf \Sigma}^x$.
\end{theorem}
{\bf Proof}\; For any increasing  twice differential function $f: \Bbb{R}^n \rightarrow \Bbb{R}$,  the ``if" part follows immediately from Corollary 2.1, since  ${\boldsymbol \mu}^x\le {\boldsymbol \mu}^y$, ${\bf \Sigma}^y={\bf \Sigma}^x$ and  $\nabla f({\bf x})\ge 0$ for all
${\bf x}\in \Bbb{R}^n$ imply that $E[f({\bf Y})]\ge E[f({\bf X})]$. To prove the  ``only if" part,  we  choose  $f$ to have the forms of
$f({\bf x})=h_1(x_i)$ and $f({\bf x})=h_2(x_i+x_j)$, where  $h_1$ and $h_2$  are any two univariate increasing functions,  it follows from  ${\bf X}\le_{st} {\bf Y}$ that
$X_i\le_{st} Y_i$ and $X_i+Y_i\le_{st} X_j+Y_j$. Note that ${\bf{X}}\sim E_n ({\boldsymbol \mu}^x,{\bf \Sigma}^x,\phi)$ and   ${\bf{Y}}\sim E_n ({\boldsymbol \mu}^y,{\bf \Sigma}^y,\phi)$ lead to  $X_i\sim E_1 ({\mu}_i^x,{\sigma}_{ii}^x,\phi)$ and $X_i+Y_i\sim E_1 ({\mu}_i^x+{\mu}_j^x,{\sigma}_{ii}^x+ {\sigma}_{jj}^x+2\sigma^x_{ij},\phi)$. By  the symmetry of elliptical distributions, all $X_i, X_j, Y_i, Y_j, X_i+Y_i$  and  $X_j+Y_j$ are supported on ${\Bbb{R}}$.  Applying Lemma 2.2 (i) we find that
${\mu}_i^x\le {\mu}_i^y$ and $\sigma^x_{ij}=\sigma^y_{ij}$ for all $1\le i,j\le n$. Hence,  ${\boldsymbol \mu}^x\le {\boldsymbol \mu}^y$ and ${\bf \Sigma}^y={\bf \Sigma}^x$. This proves the  ``only if" part.

The following result, due to  Pan et al. (2016),   generalizes  Theorem 4 in  Scarsini (1998) and  Theorem 6  in M\"uller (2001)  for  the
multivariate normal case. Here we provided a different proof.
\begin{theorem}   Let  ${\bf{X}}\sim E_n ({\boldsymbol \mu}^x,{\bf \Sigma}^x,\phi)$ and   ${\bf{Y}}\sim E_n ({\boldsymbol \mu}^y,{\bf \Sigma}^y,\phi)$. Then the following statements are equivalent:

(1)\;  ${\boldsymbol \mu}^y={\boldsymbol \mu}^x$ and ${\bf \Sigma}^y-{\bf \Sigma}^x  \succeq {\bf O}$;

(2)\; ${\bf X}\le_{cx} {\bf Y}$;

(3)\; ${\bf X}\le_{lcx} {\bf Y}$.
\end{theorem}

{\bf Proof}\; (1)$\Rightarrow$ (2). For any twice differential convex function $f: \Bbb{R}^n \rightarrow \Bbb{R}$, using Lemma 2.5 we get
$E[f({\bf Y})]\ge E[f({\bf X})]$ since $f$ is convex if and only if its Hessian matrix $H_f$ is positive semi-definite. (2)$\Rightarrow$ (3)
is obvious. (3)$\Rightarrow$ (1) is the same as  the proof of Theorem 4.1 in Pan et al. (2016).

The following result, due to  Pan et al. (2016),   generalizes Theorem 7  in M\"uller (2001)  for  the
multivariate normal case. Here we provided a different proof.
\begin{theorem}   Let  ${\bf{X}}\sim E_n ({\boldsymbol \mu}^x,{\bf \Sigma}^x,\phi)$ and   ${\bf{Y}}\sim E_n ({\boldsymbol \mu}^y,{\bf \Sigma}^y,\phi)$. Then the following statements hold:

(1)\; If ${\boldsymbol \mu}^x\le {\boldsymbol \mu}^y$ and ${\bf \Sigma}^y-{\bf \Sigma}^x   \succeq {\bf O}$, then ${\bf X}\le_{icx} {\bf Y}$.

(2)\; If ${\bf X}\le_{icx} {\bf Y}$, then ${\boldsymbol \mu}^y\ge {\boldsymbol \mu}^x$ and  ${\bf \Sigma}^y-{\bf \Sigma}^x$ is copositive, i.e., ${\bf a}'({\bf \Sigma}^y-{\bf \Sigma}^x){\bf a}\ge 0$
for all ${\bf a}\ge 0$.
\end{theorem}
{\bf Proof}\; (1). For any twice differential increasing convex function $f: \Bbb{R}^n \rightarrow \Bbb{R}$, using Lemma 2.5, together with the conditions  ${\boldsymbol \mu}^y\ge {\boldsymbol \mu}^x$ and ${\bf \Sigma}^y-{\bf \Sigma}^x   \succeq {\bf O}$, we get $E[f({\bf Y})]\ge E[f({\bf X})]$, and thus we have ${\bf X}\le_{icx} {\bf Y}$.
The proof of (2) is the same as  the proof of Theorem 4.6(2) in Pan et al. (2016).

\begin{remark} For the case of increasing convex
order, there are no   sufficient and necessary conditions   in the literature  even for  normal distributions, see M\"uller (2001) and Pan et al. (2016).
We remark that if $({\bf \Sigma}^y-{\bf \Sigma}^x){\bf z}=0$ has a positive solution, then  ${\bf a}'({\bf \Sigma}^y-{\bf \Sigma}^x){\bf a}\ge 0$
for all ${\bf a}\ge 0$  if and only if   ${\bf \Sigma}^y-{\bf \Sigma}^x \succeq {\bf O}$ (see Theorem A1).
So we get the following if-and-only-if characterization of  increasing convex
order:\\
  Assume that ${\bf{X}}\sim E_n ({\boldsymbol \mu}^x,{\bf \Sigma}^x,\phi)$,   ${\bf{Y}}\sim E_n ({\boldsymbol \mu}^y,{\bf \Sigma}^y,\phi)$
   and $({\bf \Sigma}^y-{\bf \Sigma}^x){\bf z}=0$ has a positive solution. Then ${\bf X}\le_{icx} {\bf Y}$
  if and only if  ${\boldsymbol \mu}^y\ge {\boldsymbol \mu}^x$ and ${\bf \Sigma}^y-{\bf \Sigma}^x   \succeq {\bf O}$.
\end{remark}
\begin{remark}  It is easy to see that  ${\bf \Sigma}^y-{\bf \Sigma}^x \succeq {\bf O}$ implies  ${\bf \Sigma}^y-{\bf \Sigma}^x$ is copositive, but conversely is not true.  We give an example.
Let
$${\bf \Sigma}^x =\left(\begin{array}{cc}
 \sigma^2\;& \rho_x\sigma^2 \\
 \rho_x\sigma^2\; & \sigma^2\\
\end{array}\right),\;
{\bf \Sigma}^y =\left(\begin{array}{cc}
 \sigma^2\;& \rho_y\sigma^2\\
 \rho_y\sigma^2\; & \sigma^2 \\
\end{array}\right),
$$
 where $\sigma^2>0$ and $-1\le \rho_x<\rho_y\le 1$. Then for all ${\bf a}=(a_1,a_2)'\ge 0$,  ${\bf a}'({\bf \Sigma}^y-{\bf \Sigma}^x){\bf a}=a_1a_2\sigma^2(\rho_y-\rho_x)\ge 0$.
 But for ${\bf Z}=(z_1,-z_1)'\in {\Bbb R}^2$,  ${\bf Z}'({\bf \Sigma}^y-{\bf \Sigma}^x){\bf Z}=-z_1^2\sigma^2(\rho_y-\rho_x)\le 0$.
\end{remark}

The following result    generalizes Theorem 11 in M\"uller (2001) in which   the
multivariate normal case was considered.

\begin{theorem}   Let  ${\bf{X}}\sim E_n ({\boldsymbol \mu}^x,{\bf \Sigma}^x,\phi)$ and   ${\bf{Y}}\sim E_n ({\boldsymbol \mu}^y,{\bf \Sigma}^y,\phi)$. Then the following statements are equivalent:\\
(1)\;    ${\bf X}\le_{sm} {\bf Y}$.\\
(2)\; ${\bf X}$ and ${\bf Y}$ have the same  marginal  and  $\sigma^x_{ij}\le  \sigma^y_{ij}$ for all $1\le i<j\le n$.
\end{theorem}
{\bf Proof}\; (1)$\Rightarrow$ (2).  If  ${\bf X}\le_{sm} {\bf Y}$,  then ${\bf X}$  and ${\bf Y}$ necessarily belong to the same Fr\'echet space. In particular, ${\bf X}$ and ${\bf Y}$ have the same  marginal  (see e.g. M\"uller (2000)). Since the function $f({\bf x})= x_i x_j$ is supermodular for all $1\le i<j\le n$, we see that  ${\bf X}\le_{sm} {\bf Y}$ implies $\sigma^x_{ij}\le  \sigma^y_{ij}$ for all $1\le i<j\le n$.

Since ${\bf X}\le_{sm} {\bf Y}$  if and only if $E[f({\bf X})]\le E[f({\bf Y})]$ holds for all twice differentiable
functions $f: \Bbb{R}^n \rightarrow \Bbb{R}$ satisfying $\frac{\partial^2 }{\partial x_i \partial x_j}f({\bf x})\ge 0$ for ${\bf x}\in \Bbb{R}^n$ and all $1\le i<j\le n$, the implication  (2)$\Rightarrow$ (1) follows from Lemma 2.5.

\begin{theorem}     Let  ${\bf{X}}\sim E_n ({\boldsymbol \mu}^x,{\bf \Sigma}^x,\phi)$ and   ${\bf{Y}}\sim E_n ({\boldsymbol \mu}^y,{\bf \Sigma}^y,\phi)$  be   two $n$-dimensional elliptically
distributed random vectors supported on ${\Bbb R}^n$.

(1)\; If ${\boldsymbol \mu}^x\le {\boldsymbol \mu}^y$,   $\sigma^x_{ii}=\sigma^y_{ii}$ for $i=1,2,\cdots,n$ and  $\sigma^x_{ij}\le  \sigma^y_{ij}$ for all $1\le i<j\le n$, then   ${\bf X}\le_{ism} {\bf Y}$.

(2)\; If  ${\bf X}\le_{ism} {\bf Y}$, then ${\boldsymbol \mu}^x\le {\boldsymbol \mu}^y$,   $\sigma^x_{ii}=\sigma^y_{ii}$ for $i=1,2,\cdots,n$ and  $E(X_iX_j)\le E(Y_iY_j)$ for all $1\le i<j\le n$.
\end{theorem}
{\bf Proof}\; (1).  For any twice differentiable
functions $f: \Bbb{R}^n \rightarrow \Bbb{R}$ satisfying $\frac{\partial }{\partial x_i}f({\bf x})\ge 0$ for ${\bf x}\in \Bbb{R}^n$ and all $1\le i\le n$ and $\frac{\partial^2 }{\partial x_i \partial x_j}f({\bf x})\ge 0$ for ${\bf x}\in \Bbb{R}^n$ and all $1\le i<j\le n$, using Corollary 2.1, together with the conditions  ${\boldsymbol \mu}^y\ge {\boldsymbol \mu}^x$,  $\sigma^x_{ii}=\sigma^y_{ii}$ for $i=1,2,\cdots,n$ and  $\sigma^x_{ij}\le  \sigma^y_{ij}$ for all $1\le i<j\le n$, we get $E[f({\bf Y})]\ge E[f({\bf X})]$. Thus, we have ${\bf X}\le_{ism} {\bf Y}$.

 (2).  ${\bf X}\le_{ism} {\bf Y}$ implies
$X_i\le_{st} Y_i$ (M\"uller and Stoyan (2002. P.114).    Applying Lemma 2.1 (i) we find that
${\mu}_i^x\le {\mu}_i^y$ and $\sigma^x_{ij}=\sigma^y_{ij}$ for all $1\le i,j\le n$.    Choosing a supermodular  $f({\bf x})=x_i x_j \;(i\neq j)$,   it follows from  ${\bf X}\le_{ism} {\bf Y}$ that $E(X_iX_j)\le E(Y_iY_j)$.

\begin{corollary}   Let  ${\bf{X}}\sim E_n ({\bf 0},{\bf \Sigma}^x,\phi)$ and   ${\bf{Y}}\sim E_n ({\bf 0},{\bf \Sigma}^y,\phi)$  be   two $n$-dimensional elliptically
distributed random vectors supported on ${\Bbb R}^n$.   Then the following statements are equivalent:\\
(1)\;  ${\bf X}\le_{ism} {\bf Y}$.  \\
(2)\;   $\sigma^x_{ii}=\sigma^y_{ii}$ for $i=1,2,\cdots,n$ and   $\sigma^x_{ij}\le \sigma^y_{ij}$ for all $1\le i<j\le n$.
\end{corollary}

The following result    generalizes Theorem 12 in M\"uller (2001) in which   the
multivariate normal case was considered.

\begin{theorem}   Let  ${\bf{X}}\sim E_n ({\boldsymbol \mu}^x,{\bf \Sigma}^x,\phi)$ and   ${\bf{Y}}\sim E_n ({\boldsymbol \mu}^y,{\bf \Sigma}^y,\phi)$. Then the following statements are equivalent:

(1)\;    ${\bf X}\le_{dcx} {\bf Y}$.

(2)\;  ${\boldsymbol \mu}^x={\boldsymbol \mu}^y$ and $\sigma^x_{ij}\le  \sigma^y_{ij}$ for all $1\le i,j\le n$.
\end{theorem}

{\bf Proof}\; (1)$\Rightarrow$ (2).  Note that the functions $f({\bf x})=x_i, -x_i,  x_i x_j$  are directionally convex for all $1\le i,j\le n$,
 thus ${\boldsymbol \mu}^x={\boldsymbol \mu}^y$ and $\sigma^x_{ij}\le  \sigma^y_{ij}$ for all $1\le i,j\le n$.

(2)$\Rightarrow$ (1). Since ${\bf X}\le_{dcx} {\bf Y}$  if and only if $E[f({\bf X})]\le E[f({\bf Y})]$ holds for all twice differentiable
functions $f: \Bbb{R}^n \rightarrow \Bbb{R}$ satisfying $\frac{\partial^2 }{\partial x_i \partial x_j}f({\bf x})\ge 0$ for ${\bf x}\in \Bbb{R}^n$ and all $1\le i, j\le n$, the implication  (2)$\Rightarrow$ (1) follows from Lemma 2.5.

For increasing directionally convex orders we have

\begin{theorem}   Let  ${\bf{X}}\sim E_n ({\boldsymbol \mu}^x,{\bf \Sigma}^x,\phi)$ and   ${\bf{Y}}\sim E_n ({\boldsymbol \mu}^y,{\bf \Sigma}^y,\phi)$. Then the following statements are equivalent:

(1)\;    ${\bf X}\le_{idcx} {\bf Y}$.

(2)\;  ${\boldsymbol \mu}^x\le {\boldsymbol \mu}^y$ and $\sigma^x_{ij}\le  \sigma^y_{ij}$ for all $1\le i,j\le n$.
\end{theorem}

{\bf Proof}\; (1)$\Rightarrow$ (2).  Note that the functions $f({\bf x})=x_i,   x_i x_j$  are directionally convex for all $1\le i,j\le n$.
  Thus, ${\boldsymbol \mu}^x\le {\boldsymbol \mu}^y$ and $\sigma^x_{ij}\le  \sigma^y_{ij}$ for all $1\le i,j\le n$.

(2)$\Rightarrow$ (1). Since ${\bf X}\le_{idcx} {\bf Y}$  if and only if $E[f({\bf X})]\le E[f({\bf Y})]$ holds for all twice differentiable
increasing functions $f: \Bbb{R}^n \rightarrow \Bbb{R}$ satisfying $\frac{\partial^2 }{\partial x_i \partial x_j}f({\bf x})\ge 0$ for ${\bf x}\in \Bbb{R}^n$ and all $1\le i, j\le n$, the implication  (2)$\Rightarrow$ (1) follows from Lemma 2.5.

As pointed out by  M\"uller (2001), the if-and-only-if characterization of the upper orthant order for multinormal distributions is not found.
The following result    generalizes and strengthen Theorem 10 in M\"uller (2001) in which   the
multivariate normal case was considered and Theorem 2 in  Landsman and  Tsanakas  (2006) in which   the bivariate elliptical distributions
 was considered.

\begin{theorem}   Let  ${\bf{X}}\sim E_n ({\boldsymbol \mu}^x,{\bf \Sigma}^x,\phi)$ and   ${\bf{Y}}\sim E_n ({\boldsymbol \mu}^y,{\bf \Sigma}^y,\phi)$  be   two $n$-dimensional elliptically
distributed random vectors supported on ${\Bbb R}^n$.

(1)\; If ${\boldsymbol \mu}^x\le {\boldsymbol \mu}^y$,   $\sigma^x_{ii}=\sigma^y_{ii}$ for $i=1,2,\cdots,n$ and  $\sigma^x_{ij}\le  \sigma^y_{ij}$ for all $1\le i<j\le n$, then ${\bf X}\le_{uo} {\bf Y}$.

(2)\; If  ${\bf X}\le_{uo} {\bf Y}$, then ${\boldsymbol \mu}^x\le {\boldsymbol \mu}^y$,   $\sigma^x_{ii}=\sigma^y_{ii}$ for $i=1,2,\cdots,n$ and  $E(X_iX_j)\le E(Y_iY_j)$ for all $1\le i<j\le n$.
\end{theorem}

{\bf Proof}\; (1). For any
  $\Delta$-monotone functions $f: \Bbb{R}^n \rightarrow \Bbb{R}$, using Lemma 2.5, together with the conditions  ${\boldsymbol \mu}^y\ge {\boldsymbol \mu}^x$,  $\sigma^x_{ii}=\sigma^y_{ii}$ for $i=1,2,\cdots,n$ and  $\sigma^x_{ij}\le  \sigma^y_{ij}$ for all $1\le i<j\le n$, we get $E[f({\bf Y})]\ge E[f({\bf X})]$, and thus we have ${\bf X}\le_{uo} {\bf Y}$.

 (2). Using the fact that ${\bf X}\le_{uo} {\bf Y}$ implies $X_i\le_{st} Y_i$ for all $1\le i\le n$ and Lemma 2.2(i) we get
  ${\boldsymbol \mu}^x\le {\boldsymbol \mu}^y$ and   $\sigma^x_{ii}=\sigma^y_{ii}$ for $i=1,2,\cdots,n$.
  Choosing   $f({\bf x})=x_i x_j \;(i\neq j)$, which is a    $\Delta$-monotone function,  it follows from  ${\bf X}\le_{uo} {\bf Y}$ that $E(X_iX_j)\le E(Y_iY_j)$.

%Note that $E(X_iX_j)=\sigma_{ij}^x-{\mu}_i^x {\mu}_j^x$ and
 % $E(Y_iY_j)=\sigma_{ij}^y-{ \mu}_i^y { \mu}_j^y$, if ${\boldsymbol \mu}^x ={\boldsymbol \mu}^y$,    or   ${\boldsymbol \mu}^x<{\bf 0}$ and ${\boldsymbol \mu}^y<{\bf 0}$, we get  $\sigma^x_{ij}\le  \sigma^y_{ij}$ for all $1\le i<j\le n$, since ${\mu}_i^x {\mu}_j^x\ge { \mu}_i^y { \mu}_j^y$. So we have

\begin{corollary}   Let  ${\bf{X}}\sim E_n ({\bf 0},{\bf \Sigma}^x,\phi)$ and   ${\bf{Y}}\sim E_n ({\bf 0},{\bf \Sigma}^y,\phi)$  be   two $n$-dimensional elliptically
distributed random vectors supported on ${\Bbb R}^n$. Then the following statements are equivalent:\\
(1)\;  ${\bf X}\le_{uo} {\bf Y}$.  \\
(2)\;    $\sigma^x_{ii}=\sigma^y_{ii}$ for $i=1,2,\cdots,n$ and   $\sigma^x_{ij}\le \sigma^y_{ij}$ for all $1\le i<j\le n$.
\end{corollary}

%\begin{corollary}   Let  ${\bf{X}}\sim E_n ({\boldsymbol \mu}^x,{\bf \Sigma}^x,\phi)$ and   ${\bf{Y}}\sim E_n ({\boldsymbol \mu}^y,{\bf \Sigma}^y,\phi)$  be   two $n$-dimensional elliptically
%distributed random vectors supported on ${\Bbb R}^n$. \\
%(i)   If ${\boldsymbol \mu}^x ={\boldsymbol \mu}^y$,  then
% ${\bf X}\le_{uo} {\bf Y}\Leftrightarrow {\bf X}\le_{ism} {\bf Y}$.\\
%(ii) If    ${\boldsymbol \mu}^x<{\bf 0}$ and ${\boldsymbol \mu}^y<{\bf 0}$, then
% ${\bf X}\le_{uo} {\bf Y}\Leftrightarrow {\bf X}\le_{ism} {\bf Y}$.
%\end{corollary}

The following theorem considers the componentwise convex order. The  multivariate normal case  can be found in  M\"uller and  Stoyan (2002), see also Arlotto and Scarsini (2009).

\begin{theorem}   Let  ${\bf{X}}\sim E_n ({\boldsymbol \mu}^x,{\bf \Sigma}^x,\phi)$ and   ${\bf{Y}}\sim E_n ({\boldsymbol \mu}^y,{\bf \Sigma}^y,\phi)$. Then the following statements are equivalent:

(1)\;    ${\bf X}\le_{ccx} {\bf Y}$.

(2)\;   ${\boldsymbol \mu}^x= {\boldsymbol \mu}^y$ and    $\sigma^x_{ii}\le  \sigma^y_{ii}$ for all $1\le i\le n$, and  $\sigma^x_{ij}=  \sigma^y_{ij}$ for all $1\le i<j\le n$.
\end{theorem}

{\bf Proof}\; (1)$\Rightarrow$ (2).  Note that the functions $f({\bf x})=x_i, -x_i, x_i^2, x_i x_j, - x_i x_j $  are  componentwise convex for all $1\le i,j\le n$.
 Thus, we get ${\boldsymbol \mu}^x={\boldsymbol \mu}^y$,    $\sigma^x_{ii}\le  \sigma^y_{ii}$ for all $1\le i\le n$ and $\sigma^x_{ij}=\sigma^y_{ij}$ for all $1\le i<j\le n$.

(2)$\Rightarrow$ (1).
For any twice differentiable
functions $f: \Bbb{R}^n \rightarrow \Bbb{R}$ satisfying $\frac{\partial^2 }{\partial x_i^2}f({\bf x})\ge 0$ for ${\bf x}\in \Bbb{R}^n$ and all $1\le i\le n$, using Lemma 2.5, together with the conditions   ${\boldsymbol \mu}^x= {\boldsymbol \mu}^y$,    $\sigma^x_{ii}\le  \sigma^y_{ii}$ for all $1\le i\le n$ and  $\sigma^x_{ij}=  \sigma^y_{ij}$ for all $1\le i<j\le n$,  we get $E[f({\bf X})]\le E[f({\bf Y})]$. Thus  ${\bf X}\le_{ccx} {\bf Y}$.

Similarly,  we establish the result for increasing componentwise convex order   as follows.

\begin{theorem}   Let  ${\bf{X}}\sim E_n ({\boldsymbol \mu}^x,{\bf \Sigma}^x,\phi)$ and   ${\bf{Y}}\sim E_n ({\boldsymbol \mu}^y,{\bf \Sigma}^y,\phi)$. Then the following statements are equivalent:

(1)\;    ${\bf X}\le_{iccx} {\bf Y}$.

(2)\;   $\mu_i^x\le \mu_i^y$ and    $\sigma^x_{ii}\le  \sigma^y_{ii}$ for all $1\le i\le n$, and  $\sigma^x_{ij}=  \sigma^y_{ij}$ for all $1\le i<j\le n$.
\end{theorem}
At the end of this section, we will consider the  copositive and completely positive orders for multivariate elliptical random variables.  The  multivariate normal case  can be found in   Pan et al. (2016) pointed out that it is still unknown
whether such a characterization holds for multivariate elliptical distributions. Before we state Theorem 3.11, we first give the following definitions.

{\bf Definition 3.1} (Arlotto and Scarsini (2009)) An $n\times n$ matrix ${\bf A}$ is called copositive if the quadratic form ${\bf x'Ax} \ge 0$ for all ${\bf x}\ge 0$, and ${\bf A}$ is called completely positive if there exists a nonnegative $m\times n$  matrix ${\bf B}$ such that ${\bf A}={\bf B'B}$.

Denote by ${\cal C}_{cop}$ the cone of copositive matrices and as  ${\cal C}_{cp}$ the cone of completely positive matrices. Let  ${\cal C}^*_{cop}$ and  ${\cal C}^*_{cp}$ be the dual of ${\cal C}_{cop}$ and  ${\cal C}_{cp}$, respectively. It is well known that (see Arlotto and Scarsini (2009))
  ${\cal C}^*_{cop}={\cal C}_{cp}$  and ${\cal C}^*_{cp}={\cal C}_{cop}$.

 The following Hessian orders can be defined (see Arlotto and Scarsini (2009)).\\
(a) $X \le _{cp} Y$ if   $E[f({\bf X})]\le E[f({\bf Y})]$  holds for all functions $f$ such that  ${\bf H}_f({\bf x})\in {\cal C}_{cp}$.\\
(b)  $X \le _{cop} Y$ if   $E[f({\bf X})]\le E[f({\bf Y})]$  holds for all functions $f$ such that  ${\bf H}_f({\bf x})\in {\cal C}_{cop}$.

\begin{theorem}   Let  ${\bf{X}}\sim E_n ({\boldsymbol \mu}^x,{\bf \Sigma}^x,\phi)$ and   ${\bf{Y}}\sim E_n ({\boldsymbol \mu}^y,{\bf \Sigma}^y,\phi)$. Then

(1)\;    ${\bf X}\le_{cp} {\bf Y}$ if and only if   ${\boldsymbol \mu}^x= {\boldsymbol \mu}^y$ and ${\bf \Sigma}^y-{\bf \Sigma}^x$ is copositive.

(2)\;    ${\bf X}\le_{cop} {\bf Y}$ if and only if   ${\boldsymbol \mu}^x= {\boldsymbol \mu}^y$ and ${\bf \Sigma}^y-{\bf \Sigma}^x$ is  completely copositive.
\end{theorem}
{\bf Proof}\; We prove (1)  and the proof of (2) is similar.   ``If part":  Consider the functions $f_i({\bf x})=x_i, -x_i$ $(1\le i\le n)$. Observe that  ${\bf H}_{f_i}({\bf x})={\bf O}\in {\cal C}_{cp}$.
  Thus,   ${\bf X}\le_{cp} {\bf Y}$ implies  ${\boldsymbol \mu}^x= {\boldsymbol \mu}^y$.
  Let $E({\bf X})=E({\bf Y})={\boldsymbol \mu}$.   For any    symmetric $n\times n$  matrix ${\bf A}\in {\cal C}_{cp}$, define a function $f$ as
  $$f({\bf x})=\frac12 ({\bf x}- {\boldsymbol \mu})'{\bf A}({\bf x-{\boldsymbol \mu} }).$$
  Observe that ${\bf H}_{f}({\bf x})={\bf A}$ for all ${\bf x}$, and thus ${\bf X}\le_{cp} {\bf Y}$ implies  $E[f({\bf X})]\le E[f({\bf Y})]$, which
  is equivalent to
   $$E({\bf X}- {\boldsymbol \mu})'{\bf A}({\bf X-{\boldsymbol \mu} })\le E({\bf Y}- {\boldsymbol \mu})'{\bf A}({\bf Y-{\boldsymbol \mu} }).$$
  It follows from the above that $-2\phi'(0){\rm tr}({\bf \Sigma}^x {\bf A}) \le -2\phi'(0){\rm tr}({\bf \Sigma}^y {\bf A})$.  Therefore,
  tr$(({\bf \Sigma}^y-{\bf \Sigma}^x){\bf A})\ge 0$. Since  ${\bf A}\in {\cal C}_{cp}$ is arbitrary, we conclude that  ${\bf \Sigma}^y-{\bf \Sigma}^x\in {\cal C}^*_{cp}$. Hence,
  ${\bf \Sigma}^y-{\bf \Sigma}^x$ is copositive, since ${\cal C}^*_{cp}={\cal C}_{cop}$.\\
``Only if part": For any $f$ such that ${\bf H}_f({\bf x})\in {\cal C}_{cp}$, using Lemma 2.5, together with the condition  ${\boldsymbol \mu}^x= {\boldsymbol \mu}^y$ and  the fact that ${\bf \Sigma}^y-{\bf \Sigma}^x$ is copositive, yields  $E[f({\bf X})]\le E[f({\bf Y})]$, as desired.

Concluding  the main results in this section, we have  Table 1.

\begin{table}\caption{Comparison criteria for ${\bf{X}}\sim E_n ({\boldsymbol \mu}^x,{\bf \Sigma}^x,\phi)$ and   ${\bf{Y}}\sim E_n ({\boldsymbol \mu}^y,{\bf \Sigma}^y,\phi)$}
\centering
\begin{tabular}{|c||c|c|}
\hline  Constraints on parameters & Relationship & Order\\
\hline
  ${\boldsymbol \mu}^x\le {\boldsymbol \mu}^y$, ${\bf \Sigma}^y={\bf \Sigma}^x$& $\Leftrightarrow$ &  ${\bf X}\le_{st} {\bf Y}$\\
  ${\boldsymbol \mu}^x={\boldsymbol \mu}^y$, ${\bf \Sigma}^y-{\bf \Sigma}^x  \succeq {\bf O}$& $\Leftrightarrow$& ${\bf X}\le_{cx} {\bf Y}$\\
 ${\boldsymbol \mu}^x={\boldsymbol \mu}^y$, ${\bf \Sigma}^y-{\bf \Sigma}^x  \succeq {\bf O}$ & $\Leftrightarrow$&${\bf X}\le_{lcx} {\bf Y}$\\
   ${\boldsymbol \mu}^x\le {\boldsymbol \mu}^y$, ${\bf \Sigma}^y-{\bf \Sigma}^x   \succeq {\bf O}$&  $\Rightarrow$ & ${\bf X}\le_{icx} {\bf Y}$\\
 ${\boldsymbol \mu}^x\le {\boldsymbol \mu}^y$, ${\bf \Sigma}^y-{\bf \Sigma}^x   \succeq {\bf O}$, det$({\bf \Sigma}^y-{\bf \Sigma}^x)=0$   &$\Leftrightarrow$ & ${\bf X}\le_{icx} {\bf Y}$\\
 ${\boldsymbol \mu}^x={\boldsymbol \mu}^y$, $\sigma^x_{ii}=\sigma^y_{ii}$, $\sigma^x_{ij}\le  \sigma^y_{ij}$& $\Leftrightarrow$ & ${\bf X}\le_{sm} {\bf Y}$\\
 ${\boldsymbol \mu}^x\le{\boldsymbol \mu}^y$, $\sigma^x_{ii}=\sigma^y_{ii}$, $\sigma^x_{ij}\le  \sigma^y_{ij}$& $\Rightarrow$ & ${\bf X}\le_{ism} {\bf Y}$\\
 ${\boldsymbol \mu}^x={\boldsymbol \mu}^y={\bf 0}$, $\sigma^x_{ii}=\sigma^y_{ii}$, $\sigma^x_{ij}\le  \sigma^y_{ij}$& $\Leftrightarrow$ & ${\bf X}\le_{ism} {\bf Y}$\\
 ${\boldsymbol \mu}^x={\boldsymbol \mu}^y$, $\sigma^x_{ij}\le  \sigma^y_{ij}$, $\forall i,j$ & $\Leftrightarrow $ & ${\bf X}\le_{dcx} {\bf Y}$\\
 ${\boldsymbol \mu}^x\le {\boldsymbol \mu}^y$, $\sigma^x_{ij}\le  \sigma^y_{ij}$, $\forall i,j$ & $\Leftrightarrow$ & ${\bf X}\le_{idcx} {\bf Y}$\\
 ${\boldsymbol \mu}^x\le{\boldsymbol \mu}^y$, $\sigma^x_{ii}=\sigma^y_{ii}$, $\sigma^x_{ij}\le  \sigma^y_{ij}$& $\Rightarrow$ & ${\bf X}\le_{uo} {\bf Y}$\\
  ${\boldsymbol \mu}^x={\boldsymbol \mu}^y={\bf 0}$, $\sigma^x_{ii}=\sigma^y_{ii}$, $\sigma^x_{ij}\le  \sigma^y_{ij}$& $\Leftrightarrow$ & ${\bf X}\le_{uo} {\bf Y}$\\
 ${\boldsymbol \mu}^x={\boldsymbol \mu}^y$, $\sigma^x_{ii}\le \sigma^y_{ii}$, $\sigma^x_{ij}=\sigma^y_{ij}$& $\Leftrightarrow$ & ${\bf X}\le_{ccx} {\bf Y}$\\
 ${\boldsymbol \mu}^x\le{\boldsymbol \mu}^y$, $\sigma^x_{ii}\le \sigma^y_{ii}$, $\sigma^x_{ij}=\sigma^y_{ij}$& $\Leftrightarrow$ & ${\bf X}\le_{iccx} {\bf Y}$\\
  ${\boldsymbol \mu}^x={\boldsymbol \mu}^y$, ${\bf \Sigma}^y-{\bf \Sigma}^x$ is copositive & $\Leftrightarrow$& ${\bf X}\le_{cp} {\bf Y}$\\
 ${\boldsymbol \mu}^x={\boldsymbol \mu}^y$, ${\bf \Sigma}^y-{\bf \Sigma}^x$ is completely copositive & $\Leftrightarrow$& ${\bf X}\le_{cop} {\bf Y}$\\
\hline
\end{tabular}
\end{table}

\section{ Applications and examples}

This section deals with applications of the previous results.   One
can obtain a series of probability and expectation inequalities for multivariate elliptical random variables. We
will restrict ourselves to applications concerning the supermodular ordering.

\subsection{Slepian's theorem }
Slepian's theorem for multivariate normal distributions  with non-singular covariance matrix  can be found in  Tong (1980).  Das Gupta et al. (1972)  generalized   Slepian's theorem to the elliptical distributions  with non-singular covariance matrix  which was later proved in a different way by
Joag-dev et al. (1983).  Joe (1990)  provided a shorter and elementary proof.      For its extension to  the case of singular covariance matrix the reader is referred to Fang and Liang (1989).  Here we give a simple proof.   Further results on the  normal comparison inequalities of Slepian  type can be found in Li and Shao (2002), Yan (2009) and Chernozhukov et al. (2015).

\begin{theorem}   Let  ${\bf{X}}\sim E_n ({\boldsymbol \mu}^x,{\bf \Sigma}^x,\phi)$ and   ${\bf{Y}}\sim E_n ({\boldsymbol \mu}^y,{\bf \Sigma}^y,\phi)$. If ${\bf X}$ and ${\bf Y}$ have the same marginals and  $\sigma^x_{ij}\le  \sigma^y_{ij}$ for all $1\le i<j\le n$,
then
$$P(X_1\le a_1,\cdots,X_n\le a_n)\le P(Y_1\le a_1,\cdots,Y_n\le a_n)$$
and
$$P(X_1>a_1,\cdots,X_n> a_n)\le P(Y_1> a_1,\cdots,Y_n> a_n)$$
hold for every ${\bf a}\in {\Bbb{R}}^n$. Furthermore, the inequality is strict if   $\sigma^x_{ij}\le  \sigma^y_{ij}$  for some $i, j$ and
if the supports of ${\bf X, Y}$ are unbounded.
\end{theorem}

{\bf Proof}\;  Using Theorem 3.4 and the implication ${\bf X}\le_{sm} {\bf Y}$   $\Rightarrow$ ${\bf X}\le _{uo} {\bf Y}$ and ${\bf X}\le _{lo} {\bf Y}$ yield the desired result.

The following result is an immediate consequence of Theorem 4.1.

\begin{corollary}   Let  ${\bf{X}}\sim E_n ({\bf 0},{\bf \Sigma},\phi)$,
then the probabilities
$P(\min_{1\le i\le n}g_i(X_i)>C)$ and $P(\max_{1\le i\le n}g_i(X_i)\le C)$ are increasing in each $\sigma_{ij}$, where $g_i:{\Bbb R}\rightarrow {\Bbb R}, i=1,2,\cdots, n,$ are  either all increasing or  all decreasing.
\end{corollary}

\subsection{Moment Inequalities}

In this section we  can easily derive various simple but useful inequalities for certain functions of  multivariate elliptical random variables. The proofs are based on  the results in Section  3, the most important result is the one  on  supermodular orders.
We remark that   supermodular functions play a significant role in applied fields, such as risk management, insurance, queueing,  macroeconomic
dynamics, optimization and game theory.

The following are some useful results and properties of
supermodular functions. The  proofs can be found in  B\"auerle (1997), Topkis (1998), Christofides and Vaggelatou (2004) and Marshall and Olkin (2011, P. 219).

\begin{lemma} (a)\  If $f$ is increasing and supermodular, then $\max{\{f,c\}}$ is supermodular for all $c\in {\Bbb R}$.\\
(b)\  If  $f: {\Bbb R}^n\rightarrow {\Bbb R}$ is supermodular then the function $\psi$, defined by
$\psi(x_1,x_2,\cdots, x_n)=f(g_1(x_1),\cdots,g_n(x_n)),$ is also supermodular, whenever $g_i:{\Bbb R}\rightarrow {\Bbb R}, i=1,2,\cdots, n,$ are either all increasing or  all decreasing.\\
(c)\  If $f_i$ is increasing (decreasing) on $\Bbb{R}^1$ for $i=1,2,\cdots,n$, then $f(x)=\min\{f_1(x_1),\cdots, f_n(x_n)\}=-\max\{f_1(x_1),\cdots, f_n(x_n)\}$
is supermodular on $\Bbb{R}^n$.\\
(d) $H({\bf x})=\left(\sum_{k=1}^n g_i(x_i)-t\right)^+, \left(\prod_{k=1}^n g_i(x_i)-t\right)^+$ are supermodular for any $t\ge 0$.\\
(e) If $f$ is monotonic  supermodular and $g$  is increasing and convex, then $g\circ f$ is monotonic and supermodular.\\
(f)  $H({\bf x})=  \prod_{k=1}^n \phi_i(x_i)$ is supermodular, where $\phi_i:{\Bbb R}\rightarrow {\Bbb R}, i=1,2,\cdots, n,$ are either all increasing or  all decreasing.\\
(g)   The function $f(x)=\nu(x_1+\cdots+x_n)$ is supermodular, where  $\nu$ is increasing convex.\\
(h) $H({\bf x})=-\frac{1}{n-1}\sum_{i=1}^n (X_i-\overline{X})^2$ is supermodular.\\
(i) The function $H({\bf x})=\max_{1\le k\le n}\sum_{i=1}^k X_k$ is supermodular and increasing.
\end{lemma}

 The following three theorems are immediate consequence of  Theorem 3.4 and Lemma 4.1.

\begin{theorem} Let $f$  be a convex function on $(-\infty,\infty)$. Assume that  ${\bf{X}}\sim E_n ({\boldsymbol \mu}^x,{\bf \Sigma}^x,\phi)$ and   ${\bf{Y}}\sim E_n ({\boldsymbol \mu}^y,{\bf \Sigma}^y,\phi)$. If ${\bf X}$ and ${\bf Y}$ have the same  marginal and  $\sigma^x_{ij}\le  \sigma^y_{ij}$ for all $1\le i<j\le n$, then
$$ Ef(g_1(X_1)+\cdots+g_n(X_n))\le  Ef(g_1(Y_1)+\cdots+g_n(Y_n)),$$
where $g_1, \cdots, g_n$ are monotonic in the same direction. In particular,
$$ E(X_1^3+\cdots+X_n^3)^2\le   E(Y_1^3+\cdots+Y_n^3)^2.$$
\end{theorem}

\begin{theorem}  Let  ${\bf{X}}\sim E_n ({\boldsymbol \mu}^x,{\bf \Sigma}^x,\phi)$ and   ${\bf{Y}}\sim E_n ({\boldsymbol \mu}^y,{\bf \Sigma}^y,\phi)$. If ${\bf X}$ and ${\bf Y}$ have the same marginals and  $\sigma^x_{ij}\le  \sigma^y_{ij}$ for all $1\le i<j\le n$.\\
(i)\;  Assume that $f: {\Bbb R}^n\rightarrow {\Bbb R}$ is supermodular and $g_i:{\Bbb R}\rightarrow {\Bbb R}, i=1,2,\cdots, n,$ are  either all increasing or all decreasing. Then  $ Ef(g_1(X_1),\cdots, g_n(X_n))\le  Ef(g_1(Y_1),\cdots, g_n(Y_n)).$\\
(ii)\; Assume that $f$ is increasing and supermodular. Then $E\max\{f({\bf X}),0\}\le   E\max\{f({\bf Y}),0\}$.
\end{theorem}

\begin{theorem}   Assume that  ${\bf{X}}\sim E_n ({\boldsymbol \mu}^x,{\bf \Sigma}^x,\phi)$ and   ${\bf{Y}}\sim E_n ({\boldsymbol \mu}^y,{\bf \Sigma}^y,\phi)$. If ${\bf X}$ and ${\bf Y}$ have the same  marginal and  $\sigma^x_{ij}\le  \sigma^y_{ij}$ for all $1\le i<j\le n$, then\\
(i)\;$E\prod_{k=1}^n \phi_i(X_i)\le   E\prod_{k=1}^n \phi_i(Y_i),$
where $\phi_k: {\Bbb R}\rightarrow {\Bbb R}$  are monotonic in the same
direction.\\
(ii)\; $E\min\{f_1(X_1),\cdots, f_n(X_n)\}\le E\min\{f_1(Y_1),\cdots, f_n(Y_n)\}$ and\\
$E\max\{f_1(X_1),\cdots, f_n(X_n)\}\ge E\max\{f_1(Y_1),\cdots, f_n(Y_n)\}$,
where $f_i:{\Bbb R}\rightarrow {\Bbb R}, i=1,2,\cdots, n,$ are  either all increasing or  all decreasing.\\
(iii)\; $ES_x^2\ge  ES_y^2$, where
$$S_x^2=\frac{1}{n-1}\sum_{i=1}^n (X_i-\overline{X})^2.$$
(iv)\;  If $f$ is a non-decreasing convex function, then
$$Ef\left(\max_{1\le k\le n}\sum_{i=1}^k X_k\right)\le Ef\left(\max_{1\le k\le n}\sum_{i=1}^k Y_k\right).$$
\end{theorem}

{\bf Example} 4.1 (Equicorrelated elliptical variables)
Let  ${\bf{X}}\sim E_n ({\boldsymbol \mu},{\bf \Sigma}^x,\phi)$ with ${\bf \Sigma}^x=(\sigma^x_{ij})$ such that $\sigma^x_{ii}=\sigma^2, \sigma^x_{ij}=\rho_x\sigma^2$ for $1<i<j\le n$, $\sigma^2>0, \rho_x\in [-1,1]$, and let  ${\bf{Y}}\sim E_n ({\boldsymbol \mu},{\bf \Sigma}^y,\phi)$ with ${\bf \Sigma}^y=(\sigma^y_{ij})$ such that $\sigma^y_{ii}=\sigma^2, \sigma^y_{ij}=\rho_y\sigma^2$ for $1<i<j\le n, \rho_y\in [-1,1]$.  Then
 ${\bf X}\le_{sm} {\bf Y}$ if and only if $\rho_x\le \rho_y$.
B\"auerle (1997) obtained the similar  result for normal variables and $\rho_x, \rho_y\in [0,1]$. For any supermodular function   $f: {\Bbb R}^n\rightarrow {\Bbb R}$, we deduce that the expectation $Ef({\bf X})$ is increasing in $\rho_x$. We remark that for this special correlated elliptical variable, the  supermodularity of $f$ is not  necessary. For example,    if $f: {\Bbb R}^n\rightarrow {\Bbb R}$ is twice differentiable and  satisfying $\frac{\partial^2 }{\partial x_i \partial x_j}f({\bf x})\ge 0$ for ${\bf x}\in \Bbb{R}^n$ and for some $1\le i<j\le n$,   Proposition 1 in Joag-Dev, Perlman and Pitt (1983) and its Remarks on page 454   imply that  $Ef({\bf X})$ is increasing in $\rho_x$. For example, if $\rho_x\le \rho_y$, then $E(X_1X_2 X_3^{2})\le E(Y_1Y_2 Y_3^{2})$ and $E(X_1^3X_2^3 X_3^{4})\le E(Y_1^3Y_2^3 Y_3^{4})$.

{\bf Example} 4.2 (  Serial correlated elliptical variables)
Let  ${\bf{X}}\sim E_n ({\boldsymbol \mu},{\bf \Sigma}^x,\phi)$ and  ${\bf{Y}}\sim E_n ({\boldsymbol \mu},{\bf \Sigma}^y,\phi)$ with
$${\bf \Sigma}^x =\left(\begin{array}{cccc}
 \sigma^2\;& \rho_x\sigma^2\; &\cdots\; &\rho_x^{n-1}\sigma^2 \\
 \rho_x\sigma^2\; & \sigma^2\; &\cdots\;  &\rho_x^{n-2}\sigma^2 \\
 \vdots\;&\vdots\;&\ddots\;&\vdots\;\\
 \rho_x^{n-1}\sigma^2 \; & \rho_x^{n-2}\sigma^2\; &\cdots\;  & \sigma^2\\
\end{array}\right),\;
{\bf \Sigma}^y =\left(\begin{array}{cccc}
 \sigma^2\;& \rho_y\sigma^2\; &\cdots\; &\rho_y^{n-1}\sigma^2 \\
 \rho_y\sigma^2\; & \sigma^2\; &\cdots\;  &\rho_y^{n-2}\sigma^2 \\
 \vdots\;&\vdots\;&\ddots\;&\vdots\;\\
 \rho_y^{n-1}\sigma^2 \; & \rho_y^{n-2}\sigma^2\; &\cdots\;  & \sigma^2\\
\end{array}\right),
$$
 where $\sigma^2>0$ and $\rho_x, \rho_y\in [-1,1]$. Then
 ${\bf X}\le_{sm} {\bf Y}$ if and only if $\rho_x\le \rho_y$.

{\bf Appendix 1}.

{\bf Theorem A1}\; Let $f({\bf y})={\bf y}'{\bf A}{\bf y},  {\bf y}\in \Bbb{R}^n $, where ${\bf A}$ is an $n\times n$ symmetric matrix.  If
${\bf A z}=0$ has a positive solution,
then  ${\bf y}'{\bf A}{\bf y}\ge  0$
for all ${\bf y}\ge 0$  if and only if  ${\bf y}'{\bf A}{\bf y}\ge  0$ for all ${\bf y}\in \Bbb{R}^n$.

{\bf Proof}\,  We prove the ``only if" only.  Since $f$ is quadratic, using  Taylor's expansion  we get
\begin{equation*}
f({\bf x}+t{\bf y})=f({\bf x})+t{\bf y}'\nabla f({\bf x})+t^2f({\bf y}),\;  {\bf x, y}\in \Bbb{R}^n, t \in \Bbb{R},
\end{equation*}
where $\nabla f({\bf x})=2{\bf A}{\bf x}$ is the gradient of $f$. We choose a ${\bf x}_0>0$ such that ${\bf A}{\bf x}_0=0$ and $f({\bf x}_0)=0$.   Then for any ${\bf y} \in \Bbb{R}^n$ and all sufficiently small positive $t$, one has  ${\bf x}_0+t{\bf y}\ge 0$. Thus,
$f({\bf x}_0+t{\bf y})=t^2f({\bf y})$, and
consequently $f({\bf y})\ge 0$.

{\bf Appendix 2}.

 {\bf  Proof of Lemma 2.5.}  For $0\le \lambda\le 1$, define
$$ \Psi_{\lambda}({\bf t})=\exp\left(i{\bf t}'(\lambda{\boldsymbol \mu}^y+(1-\lambda){\boldsymbol \mu}^x\right)\psi({\bf t}'(\lambda {\bf \Sigma}^y+(1-\lambda) {\bf \Sigma}^x){\bf t}), {\bf t}\in \Bbb{R}^n.$$
By using the Fourier inversion theorem
$$\phi_{\lambda}({\bf x})=\left(\frac{1}{2\pi}\right)^n\int e^{-i{\bf t}'{\bf x}}\psi_{\lambda}({\bf t})d{\bf t}.$$
The derivative of $\Psi_{\lambda}$ with respect to $\lambda$ is
\begin{eqnarray*}
\frac{\partial \Psi_{\lambda}({\bf t})}{\partial \lambda}&=&i{\bf t}'({\boldsymbol \mu}^y-{\boldsymbol \mu}^x)\exp\left(i{\bf t}'(\lambda{\boldsymbol \mu}^y+(1-\lambda){\boldsymbol \mu}^x\right)\psi({\bf t}'(\lambda {\bf \Sigma}^y+(1-\lambda) {\bf \Sigma}^x){\bf t})\\
&&+{\bf t}'({\bf \Sigma}^y- {\bf \Sigma}^x){\bf t}\exp\left(i{\bf t}'(\lambda{\boldsymbol \mu}^y+(1-\lambda){\boldsymbol \mu}^x\right)\psi'({\bf t}'(\lambda {\bf \Sigma}^y+(1-\lambda) {\bf \Sigma}^x){\bf t}),
\end{eqnarray*}
and hence,
\begin{eqnarray*}
\frac{\partial \phi_{\lambda}({\bf t})}{\partial \lambda}&=&\left(\frac{1}{2\pi}\right)^n\int e^{-i{\bf t}'{\bf x}}\frac{\partial \Psi_{\lambda}({\bf t})}{\partial \lambda} d{\bf t}\\
&=&\left(\frac{1}{2\pi}\right)^n\int e^{-i{\bf t}'{\bf x}}\Psi_{\lambda}({\bf t})i{\bf t}'({\boldsymbol \mu}^y-{\boldsymbol \mu}^x)d{\bf t}\\
&&+\left(\frac{1}{2\pi}\right)^n\int e^{-i{\bf t}'{\bf x}}{\bf t}'({\bf \Sigma}^y- {\bf \Sigma}^x){\bf t}\exp\left(i{\bf t}'(\lambda{\boldsymbol \mu}^y+(1-\lambda){\boldsymbol \mu}^x\right)\psi'({\bf t}'(\lambda {\bf \Sigma}^y+(1-\lambda) {\bf \Sigma}^x){\bf t})d{\bf t}\\
&=&-\sum_{i=1}^n(\mu_i^y-\mu_i^x)\frac{\partial \phi_{\lambda}({\bf t})}{\partial x_i}+\Delta,
\end{eqnarray*}
where
$$\Delta=\left(\frac{1}{2\pi}\right)^n\int e^{-i{\bf t}'{\bf x}}{\bf t}'({\bf \Sigma}^y- {\bf \Sigma}^x){\bf t}\exp\left(i{\bf t}'(\lambda{\boldsymbol \mu}^y+(1-\lambda){\boldsymbol \mu}^x\right)\psi'({\bf t}'(\lambda {\bf \Sigma}^y+(1-\lambda) {\bf \Sigma}^x){\bf t})d{\bf t}.$$
Note that by (2.1), there exists a random variable $R\ge 0$ such that
\begin{eqnarray*}
\psi({\bf t}'{\bf t})=E\left(E(e^{iR{\bf t}'{\bf U}^{(n)}}|R)\right)=\int_0^{\infty} {}_0F_1\left(\frac{n}{2};-\frac{r^2 ||{\bf t}||^2}{4}\right)P(R\in dr),
\end{eqnarray*}
 where
$${}_0F_1(\gamma;z)=\sum_{k=0}^{\infty}\frac{1}{(\gamma)_k}\frac{z^k}{k!}.$$
Thus, for $u>0$,
\begin{eqnarray*}
\psi'(u)&=&\int_0^{\infty} {}_0F_1\left(\frac{n}{2};-\frac{r^2 u}{4}\right)P(R\in dr)\\
&=&  \int_0^{\infty} \frac{\partial}{\partial u}{}_0F_1\left(\frac{n}{2};-\frac{r^2 u}{4}\right)P(R\in dr)\\
 &=&-\frac{1}{2n}\int_0^{\infty}{}_0F_1\left(\frac{n}{2}+1;-\frac{r^2 u}{4}\right)r^2P(R\in dr)\\
 &\equiv& -\frac{E(R^2)}{2n}\psi_1(u),
\end{eqnarray*}
where
$$\psi_1(u)=\frac{1}{E(R^2)}\int_0^{\infty}{}_0F_1\left(\frac{n}{2}+1;-\frac{r^2 u}{4}\right)r^2P(R\in dr)$$
is a characteristic generator.
Here
$$c(||{\bf t}||^2):= {}_0F_1\left(\frac{n}{2}+1;-\frac{||{\bf t}||^2 }{4}\right)$$ is
 the characteristic function of  uniform distribution in the unit sphere in   $\Bbb{R}^n$. Thus, $\Delta$  can be rewritten as
  $$\Delta=\frac{E(R^2)}{2n}\sum_{i=1}^n\sum_{j=1}^n({\sigma}^y_{ij}-{\sigma}^x_{i j})\frac{\partial^2 \phi_{1\lambda}({\bf x})}{\partial x_i \partial x_j}.$$
Define $g(\lambda)=\int_{\Bbb{R}^n} f({\bf x})\phi_{\lambda}({\bf x})d{\bf x}$, then
$E[f({\bf Y})]-E[f({\bf X})]=g(1)-g(0)=\int_0^1 g'(\lambda)d\lambda$. The result follows since
\begin{eqnarray*}
  g'(\lambda)&=&\int_{\Bbb{R}^n} f({\bf x})\frac{\partial \phi_{\lambda}({\bf x})}{\partial \lambda}d{\bf x}\\
&=&\int_{\Bbb{R}^n}({\boldsymbol \mu}^y-{\boldsymbol \mu}^x)'\nabla f({\bf x})\phi_{\lambda}({\bf x})d{\bf x}\\
&&+\frac{E(R^2)}{2n}\int_{\Bbb{R}^n}tr\{({\bf \Sigma}^{y}-{\bf \Sigma}^{x}){\bf H}_f({\bf x})\}\phi_{1\lambda}({\bf x})d{\bf x}.
\end{eqnarray*}

\noindent{\bf\Large Acknowledgements} \ The author would  like to thank Professor Xiaowen Zhou   for helpful comments on an earlier draft of the paper.
%The author are sincerely grateful for the comments from the anonymous referee, which have led to significant improvements of the present paper.
%The authors would like to thank the referees for their valuable comments and careful reading of an earlier version of this manuscript.
\noindent The research was supported by the National Natural Science Foundation of China (No. 11171179, 11571198, 11701319).

\end{document}